

\input epsf.tex


\parindent=0.7truecm
\parskip=3pt plus .5pt minus .5pt
\hsize=16truecm \vsize=22.7truecm

\let\Fi=\varphi \let\eps=\varepsilon
\font\bigbf=cmbx10 scaled \magstephalf     
\font\Bigbf=cmbx12 scaled \magstep1        
\font\letranota=cmr9


\font\TenEns=msbm10
\font\SevenEns=msbm7
\font\FiveEns=msbm5
\newfam\Ensfam
\def\Ens{\fam\Ensfam\TenEns}
\textfont\Ensfam=\TenEns
\scriptfont\Ensfam=\SevenEns
\scriptscriptfont\Ensfam=\FiveEns

\def\R{{\Ens R}}

\def\N{{\Ens N}}

\font\TenCM=cmr10
\font\SevenCM=cmr7
\font\FiveCM=cmr5
\newfam\CMfam
\def\CM{\fam\CMfam\TenCM}
\textfont\CMfam=\TenCM
\scriptfont\CMfam=\SevenCM
\scriptscriptfont\CMfam=\FiveCM


\def\unim{\,{\CM i\/}}
\def\sign{\mathop{\rm sign}\nolimits}
\def\mn{{\hbox{\fiverm min}}}
\let\lamb=\lambda
\let\Fi=\varphi
\let\eps=\varepsilon
\def\mod#1{\vert #1\vert}
\def\Mod#1{\left|#1\right|}
\def\norma#1#2{\Vert #1\Vert_{#2}}
\setbox111=\hbox to 2truemm{\hrulefill}
\setbox222=\hbox to 2truemm{\vrule height 2truemm width .4truept\hfil\vrule height 2truemm width .4truept}
\def\cqd{\vbox{\offinterlineskip\copy111\copy222\copy111}}


\newcount\numerosection
\newcount\eqnumer
\newcount\lemnumer
\newcount\fignumer
\numerosection=0
\eqnumer=0
\lemnumer=0
\fignumer=0
\def\numsection{\global\advance\numerosection by1
\global\eqnumer=0
\global\lemnumer=0
\global\fignumer=0
\the\numerosection}
\def\strutdepth{\dp\strutbox}
\def\marginalsigne#1{\strut
    \vadjust{\kern-\strutdepth\specialsigne{#1}}}
\def\specialsigne#1{\vtop to \strutdepth{
\baselineskip\strutdepth\vss\llap{#1 }\null}}
\font\margefont=cmr10 at 6pt
\newif\ifshowingMacros
\showingMacrosfalse

\def\cite#1{\csname#1\endcsname}
\def\label#1{\gdef\currentlabel{#1}}
\def\Lefteqlabel#1{\global\advance\eqnumer by 1
\label{#1}
\ifx\currentlabel\relax
\else
\expandafter\xdef
\csname\currentlabel\endcsname{(\the\numerosection.\the\eqnumer)}
\fi
\global\let\currentlabel\relax
\ifshowingMacros
\leqno\llap{%
\margefont #1\hphantom{M}}(\the\numerosection.\the\eqnumer)%
\else
\leqno(\the\numerosection.\the\eqnumer)%
\fi
}

\def\Righteqlabel#1{\global\advance\eqnumer by 1
\label{#1}
\ifx\currentlabel\relax
\else
\expandafter\xdef
\csname\currentlabel\endcsname{(\the\numerosection.\the\eqnumer)}%
\fi
\global\let\currentlabel\relax
\ifshowingMacros
\eqno(\the\numerosection.\the\eqnumer)%
\rlap{\margefont\hphantom{M}#1}
\else
\eqno(\the\numerosection.\the\eqnumer)%
\fi
}
\let\reqlabel=\Righteqlabel
\def\numlabel#1{\global\advance\eqnumer by1
\label{#1}
\ifx\currentlabel\relax
\else
\expandafter\xdef
\csname\currentlabel\endcsname{(\the\numerosection.\the\eqnumer)}
\fi
\global\let\currentlabel\relax
\ifshowingMacros
   (\the\numerosection.\the\eqnumer)\rlap{\margefont\hphantom{M}#1}
\else
(\the\numerosection.\the\eqnumer)%
\fi
}
\def\numero{\global\advance\eqnumer by1
\number\numerosection.\number\eqnumer}
\def\lemlabel#1{\global\advance\lemnumer by 1
\ifshowingMacros%
   \marginalsigne{\margefont #1}%
\else
    \relax
\fi
\label{#1}%
\ifx\currentlabel\relax%
\else
\expandafter\xdef%
\csname\currentlabel\endcsname{\the\numerosection.\the\lemnumer}%
\fi
\global\let\currentlabel\relax%
\the\numerosection.\the\lemnumer}
\def\numlem{\global\advance\lemnumer by1
\the\numerosection.\number\lemnumer}


\let\itemBibli=\item
\def\bibl#1#2\endbibl{\par{\itemBibli{#1} #2\par}}
\def\ref.#1.{{\csname#1\endcsname}}
\newcount\bib   \bib=0
\def\bibmac#1{\advance\bib by 1
\expandafter
\xdef\csname #1\endcsname{\the\bib}}
\def\BibMac#1{\advance\bib by1
\bibl{\letranota[\the\bib]}%
{\csname#1\endcsname}\endbibl}
\def\MakeBibliography#1{
\noindent{\bf #1}
\bigskip
\bib=0
\let\bibmac=\BibMac
\BiblioFil
\BiblioOrd
}

\newif\ifbouquin

\bouquinfalse
\font\tenssi=cmssi10 at 10 true pt
 at 10 true pt
\def\bibliostyle#1#2#3{{\rm #1}\ifbouquin{\tenssi #2\/}\global\bouquinfalse\else{\it #2\/}\fi{\rm #3}}


\def\BiblioFil{
\def\GSS{\bibliostyle{\letranota A.~Griffin, D.W.~Snoke and Stringari (ed), Bose-Einstein Condensation, Cambridge University Press, (1995).}{}{}}
\def\EinsteinA{\bibliostyle{\letranota A.~Einstein, Sitzungsber.~K.~Preuss.~Akad.~Wiss.~Phys.~Math., { 261}, (1924).}{}{}}
\def\EinsteinB{\bibliostyle{\letranota A.~Einstein, Sitzungsber.~K.~Preuss.~Akad.~Wiss.~Phys.~Math., { 3}, (1925).}{}{}}
\def\AEMW{\bibliostyle{\letranota M.H.J.~Anderson, J.R.~Ensher, M.R.~Matthews and C.E.~Wiemna, Science, { 269}, (1995), 198.}{}{}}
\def\DMADDKK{\bibliostyle{\letranota K.B.~Davis at al., Phys.~Rev.~Lett., { 75}, (1995), 3969.}{}{}}
\def\Gross{\bibliostyle{\letranota E.P.~Gross, J.~Math.~Phys., { 4}, (1963), 195.}{}{}}
\def\Pitaev{\bibliostyle{\letranota L.K.~Pitaevskii, Sov.~Phys.~JETP, { 13}, (1961), 451.}{}{}}
\def\CipoKavian{\bibliostyle{\letranota R.~Cipolatti and O.~Kavian, J.~Diff.~Equations, { 176}, (2001), 223.}{}{}}
\def\Weinstein{\bibliostyle{\letranota M.I.~Weinstein, Comm.~Math.~Phys., { 87}, (1983), 567.}{}{}}
\def\Carles{\bibliostyle{\letranota R.~Carles, Ann.~Henri Poincar\'e, { 3}, (2003), 757.}{}{}}
\def\Zhang{\bibliostyle{\letranota J.~Zhang, Z.~angew.~Math.~Phys, {51}, (2000), 498.}{}{}}
\def\CazeLions{\bibliostyle{\letranota T.~Cazenave T and P-L.~Lions, Comm.~Math.~Phys., { 85}, (1982), 153.}{}{}}
\def\TrallA{\bibliostyle{\letranota J.D.~Drake-P\'erez at al., Phys.~Stat.~Sol.~(C), { 2}, No.~10, (2005), 3665.}{}{}}
\def\TrallB{\bibliostyle{\letranota C.~Trallero-Giner, J.C.~Drake-P\'erez, V.~L\'opez-Richard and J.L.~Birman,  Physica D, { 237}, (2008), 2342.}{}{}}
\def\TrallD{\bibliostyle{\letranota C.~Trallero-Giner, V.~L\'opez-Richard, M-C.~Chung and A.~Buchleitner, Physical Review A, { 79}, (2009), 063621.}{}{}}
\def\LiebEtAl{\bibliostyle{\letranota E.H.~Lieb, R.~Seiringer and J.~Yngvason, Physical Review A, { 61}, (2000), 043602.}{}{}}
\def\LiebLoss{\bibliostyle{\letranota E.H.~Lieb and M.~Loss, Analysis 2nd Edition, Graduate Studies in Mathematics, Vol.~14, AMS, (2001).}{}{}}
\def\Mossino{\bibliostyle{\letranota J.~Mossino, Inegalit\'es Isoperimetriques et Aplications en Physique, Hermann, (1984).}{}{}}
\def\Carretero{\bibliostyle{\letranota R.~Carretero-Gonz\'alez, D.J.~Frantzeskakis and P.G.~Kevrekidis, Nonlinearity, { 21}, (2008), R139.}{}{}}
\def\Oliver{\bibliostyle{\letranota O.~Morsch and M.~Oberthaler, Reviews of Modern Physics, Vol.~78, (2006), 179.}{}{}}
\def\Burger{\bibliostyle{\letranota S.~Burger et al., Phys.~Rev.~Lett., { 86}, (2001), 4447.}{}{}}
\def\KavWeis{\bibliostyle{\letranota O.~Kavian and F.~Weissler, Mich.~J.~Math., 41, (1994), 151.}{}{}}
\def\Oh{\bibliostyle{\letranota Y-G.~Oh, J.~Diff.~Eq., {81}, (1989), 255.}{}{}}
}

\def\BiblioOrd{
\bibmac{GSS}       
\bibmac{EinsteinA} 
\bibmac{EinsteinB} 
\bibmac{AEMW}      
\bibmac{DMADDKK}   
\bibmac{Oliver}    
\bibmac{Gross}     
\bibmac{Pitaev}    
\bibmac{LiebLoss}  
\bibmac{Zhang}     
\bibmac{Carretero} 
\bibmac{TrallB}    
\bibmac{TrallD}    
\bibmac{CipoKavian}
\bibmac{KavWeis}   
\bibmac{LiebEtAl}  
\bibmac{Mossino}   
\bibmac{Carles}    
\bibmac{Oh}        
\bibmac{CazeLions} 
\bibmac{TrallA}    
\bibmac{Burger}    

}

\BiblioOrd




\centerline{\Bigbf Bose-Einstein condensates in optical lattices: mathematical }
\centerline{\Bigbf analysis and analytical approximate formulae}

\bigskip
\centerline{R.~Cipolatti$^{1}$\footnote{*}{\letranota Corresponding author; E-mail address: cipolatti@im.ufrj.br.},  J.~L\'opez Gondar$^{1}$ and C.~Trallero-Giner$^{2}$ }
\bigskip

{\baselineskip=14truept
\letranota
\centerline{$^{1}\,$Instituto de Matem\'atica, Universidade Federal do Rio de Janeiro}
\centerline{C.P. 68530, Rio de Janeiro, RJ, Brasil}
\centerline{$^{2}\,$Faculty of Physics, Havana University}
\centerline{10400 Havana, Cuba}

}


\bigskip
\hrule
\bigskip

{
\letranota\baselineskip=12truept
\noindent{\bf Abstract}

We show that the GPE with cubic nonlinearity, as a model to describe the one dimensional Bose-Einstein condensates loaded into a harmonically
confined optical lattice, presents a set of ground states which is orbitally stable for any value of the self-interaction (attractive and repulsive) 
parameter and laser intensity. We also derive a new formalism which gives explicit expressions for the minimum energy $E_\mn$ and the associated 
chemical potential  $\mu_0$. Based on these formulas, we generalize the variational method to obtain approximate solutions, at any order of 
approximation, for  $E_\mn$, $\mu_0$ and the ground state.

\medskip
\noindent{\letranota Key words: Bose-Einstein condensates, stability of ground states, analytical approximate formulae, repulsive or attractive interatomic interactions.}

}
\bigskip
\hrule
\bigskip

\noindent{\bigbf\numsection.\ Introduction}\par
\bigskip

\noindent The Bose-Einstein condensation (BEC) is a fundamental phenomenon connected to superfluidity in liquid helium
[\cite{GSS}]. Its achievement in practice leads to the first unambiguous manifestation of the existence of a macroscopic quantum
state in a many-particle system.  Although the condensates of boson particles was firstly predicted by
Einstein [\cite{EinsteinA},\cite{EinsteinB}] in 1924,  BECs were experimentally realized only in 1995 [\cite{AEMW},\cite{DMADDKK}].
For these reasons, this phenomenon has attracted the attention of many scientists, in particular, during recent years.
Nowadays, one of the most interesting problems in cold matter physics is the study of the BEC in a potential trap loaded in a
periodical optical lattice (for a detailed discussion see Ref.~[\cite{Oliver}]).
This problem can be described by means of the Gross-Pitaevskii equation (GPE) [\cite{Gross},\cite{Pitaev},\cite{LiebLoss}], for which
the BEC are characterized by its ground state solutions.

The ground state solutions of the GPE in an external potential are not necessarily stable, i.e., small initial
perturbations around a ground state could give rise to solutions collapsing in finite time. In order to understand the properties
of the BEC, it is important to know how the small-amplitude excitations of the ground state evolve in time. One step in this direction
was given by Zhang [\cite{Zhang}], who proved the stability of these solutions in the case of a harmonic potential for a
nonlinear attractive interaction.

Assuming that the harmonic trapping potential has a strong aniso\-tropy (of ``cigar-shaped'' type),
the 1D limit of the GPE with cubic nonlinearity can be considered as a model to describe the condensate, more precisely, by the equation (see  [\cite{Carretero}])
$$-{\hbar^2\over 2m}{d^2\Psi \over d x^2}+{1\over 2}m\omega^2x^2\Psi+
    \lambda_{1D}|\Psi|^2\Psi-V_L\cos^2\left({2\pi\over d}x\right)\Psi=\mu_0\Psi,\reqlabel{GPE}$$
where $\omega>0$ is the oscilator trap frequency, $m>0$ is the atomic mass, $V_L>0$ is the laser intensity, $d>0$ is the wavelenght of the laser,
$\mu_0\in\R$ is the chemical potential and $\lambda_{1D}\in\R$ is the self-interaction
parameter ($\lambda_{1D}<0$ when interatomic forces are attractive and  $\lambda_{1D}>0$ for repulsive interatomic forces).
In its dimensionless form the equation \cite{GPE} can be written as
$$-{d^2\psi\over d\xi^2}+\xi^2\psi+\lamb\mod{\psi}^2\psi-V_0\cos^2(\alpha\xi)\psi=\mu\psi,\reqlabel{GPE1}$$
where, for $l:=\sqrt{\hbar/m\omega}$, we set $\xi={x/l}$, $\psi(\xi):=\sqrt{l}\Psi(x)$ and
$$\mu:={2\mu_0\over \hbar\omega},\,\,\lamb:={2\lambda_{1D}\over l\hbar\omega},
    \,\,V_0:={2V_L\over \hbar\omega},\,\, \alpha:={2\pi l\over d}.$$
The solutions of $\cite{GPE1}$  can be viewed as standing waves of the
time dependent GPE, namely,
$$\unim{\partial{u}\over \partial\tau}  = -{\partial^2{u}\over \partial\xi^2}
   +\xi^2u+\lamb|u|^2{u}-V_0\cos^2(\alpha\xi)u,\reqlabel{EqEvolution}$$
where $\tau:=\omega t/2$ and $t$ is the time.
By standing waves we mean time periodic solutions  of the form
$$u(\tau,\xi): =e^{-\unim\mu\tau}\psi(\xi).\reqlabel{DefStanWav}$$

In this paper we prove the ground state stability (beyond the Bogoliubov approximation) concerning the solutions of \cite{EqEvolution},
in both cases, attractive and repulsive, and for any value of $V_{0}$. This result, together with some qualitative properties
of the ground states and the general behavior of the minimal energy as function of $\lamb$, completes the essential of Section~2.
Furthermore, taking into account the importance from the physical standpoint
to have explicit formulae for the minimal energy, the corresponding chemical potential
and the ground state, we develop in Section~3 a new formalism to obtain explicit approximate expressions for the above
mentioned magnitudes. The present formalism leads to a new and more general variational approach. To verify the validity of
the method, we compare our solutions with those numerical results reported in [\cite{TrallB}, \cite{TrallD}].


\bigskip
\noindent{\bf\numsection.\ Existence and stability of ground states}\par
\bigskip

Although the solutions of \cite{EqEvolution} are in general complex valued functions, we can restrict our analysis
of the existence of ground states only for real valued ones, as we can see by the following lemma, where $H^1(\R)$ denotes the
usual Sobolev space.

\smallskip
\noindent{\bf Lemma\ \lemlabel{Lemma1}:} {\sl If $\psi\in H^1(\R)$ is a complex solution of \cite{GPE1}, then there exists a real function
$U(\xi)$ which is a solution of $\cite{GPE1}$ and a real number $\theta$ such that $\psi(\xi)=e^{\unim\theta}U(\xi)$.}

\noindent
{\bf Proof:} Assume that $\psi$ is a complex solution of \cite{GPE1} and consider
its real and imaginary parts, i.e., $\psi=u+\unim v$,  $u\not=0$. Then, it follows that
$$\eqalign{
-{d^2 u \over d\xi^2}+\xi^2 u+
    \lamb (u^2+v^2)u-V_0\cos^2\left(\alpha\xi\right)u & =\mu u,\cr
-{d^2 v \over d\xi^2}+\xi^2 v+
    \lamb (u^2+v^2)v-V_0\cos^2\left(\alpha\xi\right)v & =\mu
    v.\cr
}\reqlabel{SistRealSol}$$
Now, multiplying the first equation in the above system by $v$, the second by $u$ and subtracting, we get
$$u{d^2v\over d\xi ^2}-v{d^2u\over d\xi ^2}=0\quad\Rightarrow\quad u{dv\over d\xi }-v{du\over d\xi }=C,$$
for some $C\in\R$. Since the solutions of \cite{GPE1} tends to zero as $\xi\rightarrow\pm\infty$ (see the next Theorem),
we conclude that $C=0$ and hence
$$u{dv\over d\xi }-v{du\over d\xi }=0\quad\Rightarrow\quad {d\hfil\over d\xi }\left({v\over u}\right)=0.$$
Therefore,  $v=\gamma u$, $\gamma\not=0$ and each one of the equations of \cite{SistRealSol} reduces to
$$-{d^2 u \over d\xi^2}+\xi^2 u+
    \lamb (1+\gamma^2)u^3-V_0\cos^2\left(\alpha\xi\right)u =\mu u.$$
Now, considering $U(\xi):=(1+\gamma^2)^{1/2}u(\xi)$, it follows that $U(\xi)$ is a real valued solution of \cite{GPE1} and
$$\psi=u+\unim v=\left({1\over\sqrt{1+\gamma^2}}+{\unim\gamma\over\sqrt{1+\gamma^2}}\right)U=e^{\unim\theta}U,$$
where $\theta:=\arctan\gamma$.\quad\cqd

An important property of the solutions of \cite{GPE1} is their asymptotic decay at infinity, as asserted in the following result.

\noindent {\bf Theorem\ \lemlabel{ExpDecay}:} {\sl Let $\mu, \lamb, V_0\in\R$ be given. If $\varphi\in H^1(\R)$ is a solution of \cite{GPE1},
then $\varphi\in C^2(\R)$ and for $\delta\in(0,1)$ there exists $C(\delta)>0$ such that
$$\forall\xi\in\R,\quad \mod{\varphi(\xi)}\le C(\delta)\exp[-(1-\delta)\xi^2/2].\reqlabel{DecaiExp}$$
Moreover, if $\lamb>0$ and $\mu<1-\mod{V_0}$, the above inequality holds for $\delta=0$.}

\noindent
{\bf Proof:} We proceed as in [\cite{CipoKavian}]. Since $\varphi\in H^1(\R)$, we have that $\varphi(\xi)$ is a continuous
function satisfying
$$\lim_{\mod{\xi}\to+\infty}\varphi(\xi)=0\reqlabel{Tend2Zero}$$
and it follows directly from the equation that $\psi''\in C(\R)$. In order to prove the exponential decay of $\varphi$, let
$$a(\xi):=\xi^2-\mu+\lamb\mod{\varphi(\xi)}^2-V_0\cos^2(\alpha\xi).$$
By Kato's inequality, if $z(\xi):=\mod{\varphi(\xi)}$, we have $z''\ge\sign(\varphi)\varphi''$ in the sense of distributions.
Therefore, $-z''+a(\xi)z\le 0$ in the same sense.
On the other hand, if we set $\psi_0(\xi):=C\exp[-(1-\delta)\xi^2/2]$ for $0\le \delta<1$ and $C>0$, a simple calculation gives
$$-\psi_0''+a(\xi)\psi_0=\bigl[\delta(2-\delta)\xi^2-\mu+1-\delta+\lamb\mod{\varphi(\xi)}^2-V_0\cos^2(\alpha\xi)\bigr]\psi_0.$$
If $\lamb>0$ and $\mu<1-\mod{V_0}$ set $\delta=0$, otherwise assume that $0<\delta<1$. Then, for $R>0$ large enough, it follows from \cite{Tend2Zero}
that $a(\xi)\ge 1$ and
$$\delta(2-\delta)\xi^2-\mu+1-\delta+\lamb\mod{\varphi(\xi)}^2-V_0\cos^2(\alpha\xi)>0,\quad\forall \mod{\xi}>R,$$
so that $-z''+az\le -\psi_0''+a\psi_0$ for $\mod{\xi}>R$. Moreover, if we choose $C>0$ such that $z(\pm R)\le \psi_0(\pm R)$,
then from the maximum principle we infer that
$$\mod{\varphi(\xi)}=z(\xi)\le \psi_0(\xi),\quad\forall\mod{\xi}\ge R,$$
which implies the exponential decay of $\varphi$, as asserted in \cite{DecaiExp}.\quad\cqd

\smallskip
\noindent$\bullet$ $\underline{\hbox{\sl Existence of ground states}}$
\smallskip

We introduce the variational problem which allows to prove the existence and stability of ground states for Eq.~\cite{GPE1}. Let
$$X:=\Bigl\{\psi\in H^1(\R)\,;\, \int_\R\left(\Mod{\psi'(\xi)}^2+ \xi^2\mod{\psi(\xi)}^2\right)d x<+\infty\Bigr\}.\reqlabel{DefX}$$
$X$ is a real Hilbert space if endowed with the following
usual inner product
$$(\phi|\psi)_X:=\int_\R\left(\psi'(\xi)\phi'(\xi)+\xi^2\psi(\xi)\phi(\xi)\right)\,d\xi.$$
Then, the associated norm is given by
$$\norma{\psi}{X}^2:=\int_\R\left(\Mod{\psi'(\xi)}^2+\xi^2\mod{\psi(\xi)}^2\right)\,d\xi .$$
Now we define the ``energy'' $E:X\rightarrow\R$ and the ``charge''
$Q:X\rightarrow\R$ respectively by
$$\eqalign{
E(\psi) & := \int_\R\mod{\psi'(\xi)}^2d\xi
   +\int_\R \xi^2\mod{\psi(\xi)}^2d\xi +{\lamb \over 2}\int_\R\mod{\psi(\xi)}^4d\xi
       -V_0\int_\R\cos^2\left(\alpha\xi\right)\mod{\psi(\xi)}^2d\xi ,\cr
Q(\psi) & := \int_\R\mod{\psi(\xi)}^2d\xi ,\cr }\reqlabel{DefEQ}$$
and the manifold
$$\Sigma_1:=\Bigl\{\psi\in X\,;\, Q(\psi)=1\Bigr\}.$$
With these ingredients we look for solutions  $\psi$ of Eq.~\cite{GPE1} that minimizes the energy $E$ among all functions in $\Sigma_1$.
More precisely, we look for $\psi_\mn\in\Sigma_1$ such that
$$E(\psi_\mn)=\min\bigl\{E(\psi)\,;\, \psi\in\Sigma_1\bigr\}.\reqlabel{VariatProb1}$$

\noindent{\bf Remark\ \lemlabel{Rem3}:} Before proceeding to prove that there exist
solutions of the variational problem \cite{VariatProb1}, we
shall remember the following well known facts. Let $\varphi_0:\R\rightarrow\R$ be the function
$$\varphi_0(\xi):={1\over \root 4 \of\pi}\exp(-\xi^2/2).\reqlabel{PsiZero}$$
It is easy to see that $\varphi_0\in\Sigma_1$ and that $-\varphi_0''(\xi)+\xi^2\varphi_0(\xi)=\varphi_0(\xi)$. This means that $\varphi_0$ is
an eigenfunction of the operator $L=-{d^2\hfil\over d\xi ^2}+\xi^2$ corresponding to the eigenvalue $\lambda_0=1$. In fact, $L$ has an
infinite sequence of eigenvalues $\lambda_0<\lambda_1<\cdots$, where $\lambda_n=(2n+1)$, ($n=0,1,\ldots$), and the
Hermite functions are the corresponding eigenfunctions. It is also known that $\lambda_0$ has the following variational characterization,
$$\lambda_0=\inf\Bigl\{\int_\R\left(\Mod{\psi'(\xi)}^2+\xi^2\mod{\psi(\xi)}^2\right)\,d\xi \,;\, \psi\in\Sigma_1\Bigr\}$$
and we can easily verify that $\lambda_0=1$ and that the above infimum is actually a minimum attained at $\varphi_0$, i.e.,
$\lambda_0=\norma{\varphi_0}{X}^2=1$.

We are now in position to prove the existence of ground states of \cite{GPE1}.

\noindent
{\bf Theorem\ \lemlabel{Thm1}:} {\sl Let $\lamb, V_0\in\R$ be given. Then, there exists $\psi_\mn\in\Sigma_1$ such that }
$$E(\psi_\mn)=\min\Bigl\{E(\psi)\,;\, \psi\in\Sigma_1\Bigr\}.\reqlabel{Variational}$$

\smallskip\noindent
{\bf Proof:} We divide the proof in three steps.

\noindent{\bf Step 1:} {\sl The energy $E$ is bounded from bellow on $\Sigma_1$\/}:

From  Remark~\cite{Rem3} we know that $\norma{\psi}{X}^2\ge 1$ for all $\psi\in\Sigma_1$. Hence, if $\lamb \ge0$, we have
$$E(\psi)\ge \norma{\psi}{X}^2-\mod{V_0}\ge 1-\mod{V_0},\quad\forall \psi\in\Sigma_1.\reqlabel{PrimaLimit}$$
 On the other hand, from Gagliardo-Nirenberg inequality, there
exists a constant $C_{\hbox{\letranota gn}}>0$ such that
$$\norma{\psi}{4}^4\le C_{\hbox{\letranota gn}}\norma{\psi}{2}^3\norma{\psi'}{2},\quad\forall\psi\in H^1(\R),\reqlabel{Gal-Nir}$$
where $\norma{\cdot}{4}$ and $\norma{\cdot}{2}$ are the standard norms of the spaces $L^4(\R)$ and $L^2(\R)$, respectively.
Hence, in the attractive case ($\lamb <0$), we have for any $\psi\in\Sigma_1$,
$$ E(\psi) \ge \norma{\psi}{X}^2+{\lamb C_{\hbox{\letranota gn}}\over 2}\norma{\psi'}{2}-\mod{V_0}
              \ge \norma{\psi}{X}^2+{\lamb C_{\hbox{\letranota gn}}\over 2}\norma{\psi}{X}-\mod{V_0},\reqlabel{BigDesig}$$
from which we get
$$E(\psi)\ge -{\lamb^2 C_{\hbox{\letranota gn}}^2\over 16}-\mod{V_0},\quad\forall\psi\in\Sigma_1.\reqlabel{DuaLimit} $$

\smallskip
\noindent{\bf Step 2:} {\sl The variational problem \cite{Variational} has a solution\/}:

Let $E_\mn:=\inf\bigl\{E(\psi)\,;\, \psi\in\Sigma_1\bigr\}$. From Step~1, it follows that $E_\mn\in\R$ and from the definition of infimum, we
conclude that there exists a sequence of minimizing functions $\{\psi_n\}_{n\in\N}$ in the manifold $\Sigma_1$, i.e.,
$$\forall n\in\N,\,\,\exists \psi_n\in\Sigma_1\quad\hbox{\rm such that}\quad\lim_{n\rightarrow+\infty}E(\psi_n)=E_\mn.$$
Assuming that $\lamb \ge 0$, we obtain easily from \cite{PrimaLimit} that $\{\psi_n\}$
is a bounded sequence in $X$. On the other hand, if $\lamb <0$, from \cite{BigDesig} and the Young's inequality, we have
$$\norma{\psi_n}{X}^2  \le E(\psi_n)-{\lamb  C_{\hbox{\letranota gn}}\over 2}\norma{\psi_n}{X}+\mod{V_0}
                            \le E(\psi_n)+{1\over 2}\norma{\psi_n}{X}^2+{\lamb ^2C_{\hbox{\letranota gn}}^2\over 8}+\mod{V_0}.$$
So, we obtain
$${1\over 2}\norma{\psi_n}{X}^2\le E(\psi_n)+{\lamb^2C_{\hbox{\letranota gn}}^2\over 8}+\mod{V_0},$$
from which we conclude that $\{\psi_n\}$ is a bounded sequence in $X$.
Therefore, in both cases, it follows from the Banach-Alaoglu Theorem, that there exists a subsequence of $\psi_{n}$ that converges to some
$\psi_\mn$ in the weak topology of $X$, i.e., $\psi_{n_k}\rightharpoonup\psi_\mn$. To simplify the notation, we still write $\psi_n$ for this
subsequence.
Since the embedding $X\subset L^p(\R)$ (see [\cite{CipoKavian},\cite{KavWeis}]) is compact for all $2\le p<+\infty$, we have
$$\eqalign{
\lim_{n\rightarrow+\infty}\int_\R\mod{\psi_n(\xi)}^4d\xi  & = \int_\R\mod{\psi_\mn(\xi)}^4d\xi ,  \cr
\lim_{n\rightarrow+\infty}\int_\R\mod{\psi_n(\xi)}^2d\xi  & = \int_\R\mod{\psi_\mn(\xi)}^2d\xi ,  \cr
\lim_{n\rightarrow+\infty}\int_\R\cos^2\left(\alpha\xi\right)\mod{\psi_n(\xi)}^2d\xi  & =
                          \int_\R\cos^2\left(\alpha\xi\right)\mod{\psi_\mn(\xi)}^2d\xi. \cr
}\reqlabel{Part1}$$
By hypothesis, $\psi_n\in\Sigma_1$ for all $n\in\N$, and the second limit above implies that $\psi_\mn\in\Sigma_1$.
Moreover, as the norm $\norma{\ }{X}$ is semi-continuous for the weak topology of $X$, we have
$$\norma{\psi}{X}^2\le  \liminf_{n\rightarrow+\infty}\norma{\psi_n}{X}^2.\reqlabel{Part2}$$
From \cite{Part1} and \cite{Part2}, we conclude that $E(\psi_\mn)=E_\mn$ and hence $\psi_\mn$ is a solution of \cite{Variational}.

\smallskip
\noindent {\bf Step 3:}   {\sl The function $\psi_\mn$ is a solution of \cite{GPE1}. }

This follows directly from the fact that $E(\psi)$ and $Q(\psi)$ are differentiable functionals in $X$. Indeed, from the Lagrange Theorem,
there exists $\mu\in\R$ (a Lagrange multiplier) such that
$$E'(\psi_\mn)=\mu Q'(\psi_\mn),$$
where $E'(\psi)$ and  $Q'(\psi)$ are the Fr\'echet derivatives of $E$ and $Q$ at $\psi$, respectively.
Note that this last equation is the same as \cite{GPE1}. This completes the proof. \quad\cqd

\smallskip
If we denote by ${\cal G}$ the set of ground states of \cite{GPE1}, i.e.,
$${\cal G}:=\Bigl\{\psi\in\Sigma_1\,;\, E(\psi)=E_\mn\Bigr\},$$
it follows easily from \cite{PrimaLimit} and \cite{BigDesig} that $\cal G$ is a bounded set of $X$ and
we have the following properties:

\noindent
{\bf Theorem\ \lemlabel{Thm1-1}:} {\sl Let $\lamb, V_0\in\R$.
{\parindent=0.7truecm
\item{a)} If $\lamb\ge 0$, there exists a unique positive symmetric function $\psi_\mn\in\Sigma_1$ such that
$${\cal G}=\bigl\{e^{\unim\theta}\psi_\mn\,;\, \theta\in\R\bigr\}.$$
\item{b)} If $V_0=0$, then there exists a positive symmetric function $\psi_\mn\in{\cal G}$ such that
 $\xi\mapsto\psi_\mn(\xi)$ is decreasing in the interval $\xi\ge 0$.
In particular, $\psi_\mn(0)=\max\{\psi_\mn(\xi)\,;\, \xi\in\R\}$.

}}

\smallskip\goodbreak\noindent
{\bf Proof:} To prove (a) we proceed as in [\cite{LiebEtAl}]. Suppose that there are two real functions $\psi_0,\psi_1\in{\cal G}$,
 $\psi_0\not=\psi_1$. If $\mod{\psi_0}\not=\mod{\psi_1}$, define $\psi_\nu\colon=[\nu\psi_1^2+(1-\nu)\psi_0^2]^{1/2}$, where $0<\nu<1$.
 It follows that $\psi_\nu\in\Sigma_1$ and
$$\eqalign{
\int_\R \xi^2\mod{\psi_\nu(\xi)}^2\,d\xi & = \nu\int_\R \xi^2\mod{\psi_1(\xi)}^2\,d\xi
          +(1-\nu)\int_\R \xi^2\mod{\psi_0(\xi)}^2\,d\xi,\cr
\int_\R \cos^2(\alpha\xi)\mod{\psi_\nu(\xi)}^2\,d\xi & = \nu\int_\R \cos^2(\alpha\xi)\mod{\psi_0(\xi)}^2\,d\xi
          +(1-\nu)\int_\R \cos^2(\alpha\xi)\mod{\psi_1(\xi)}^2\,d\xi.\cr
}\reqlabel{Form*1}$$
Since $x\mapsto \mod{x}^4$ is strictly convex, we have
$$\int_\R \mod{\psi_\nu(\xi)}^4\,d\xi  < \nu\int_\R \mod{\psi_0(\xi)}^4\,d\xi
          +(1-\nu)\int_\R \mod{\psi_1(\xi)}^4\,d\xi.\reqlabel{Form*2}$$
Moreover, by differentiating both sides of $\psi_\nu^2=\nu\psi_0^2+(1-\nu)\psi_1^2$ and using
the Cauchy-Schwarz inequality in $\R^2$ ($ab+cd\le \sqrt{a^2+c^2}\sqrt{b^2+d^2}$), we get
$$\psi_\nu\psi_\nu'  =   \nu\psi_0\psi_0' + (1-\nu)\psi_1\psi_1' \le \psi_\nu\sqrt{\nu(\psi_0')^2+(1-\nu)(\psi_1')^2},$$
from which it follows that
$$\mod{\psi_\nu'(\xi)}^2\le \nu\mod{\psi_0'(\xi)}^2+(1-\nu)\mod{\psi_1'(\xi)}^2.\reqlabel{Form*3}$$
As we are assuming that $\lamb\ge 0$, it follows from \cite{Form*1}--\cite{Form*3} that $E(\psi_\nu)<\nu E_\mn+(1-\nu)E_\mn=E_\mn$, which is impossible by the
definition of $E_\mn$.

Since $\mod{\psi_i}\in\Sigma_1$ for $i=0,1$, it follows from Kato's inequality that $E(\mod{\psi_i})\le E(\psi_i)=E_\mn$.
So, $\mod{\psi_0},\mod{\psi_1}\in{\cal G}$ and we conclude by the previous arguments that $\mod{\psi_0}=\mod{\psi_1}$.
By assuming that $\psi_0(\xi_0)=\psi_1(\xi_0)=0$ for some $\xi_0\in\R$ and taking into account that $\mod{\psi_i}\in C^2$ (see Theorem~\cite{ExpDecay}),
it follows that $\psi_0'(\xi_0)=\psi_1'(\xi_0)=0$, which implies from the uniqueness of solutions of ODEs that $\psi_i\equiv 0$.
This is in  contradiction from the fact that $\psi_i\in\Sigma_1$. Therefore, we can assume that $\psi_1=-\psi_0$, with $\psi_0(\xi)> 0$ for all $\xi\in\R$.
Moreover, since $\widetilde\psi(\xi)\colon=\psi(-\xi)$ belongs to $\cal G$ for all $\psi\in\cal{G}$, the unique positive ground state $\psi_0$ is
necessarily symmetric.

To prove (b), let $\psi\in{\cal G}$ be a real function and consider $\psi_*(\xi)$ the {\sl symmetric-decreasing rearrangement\/}
of $\mod{\psi(\xi)}$. It is well known (see [\cite{LiebLoss},\cite{Mossino}]) that $\psi_*$ is positive, symmetric, decreasing
in $[0,+\infty)$ and
$$\eqalign{
\int_\R \mod{\psi_*(\xi)}^p\,d\xi  & = \int_\R \mod{\psi(\xi)}^p\,d\xi,\,\, 1\le p\le\infty, \cr
\int_\R \mod{\psi'_*(\xi)}^2\,d\xi & \le \int_\R \mod{\psi'(\xi)}^2\,d\xi. \cr
}\reqlabel{Form*4}$$
Hence, from the first equality of \cite{Form*4} with $p=2$ we get $\psi_*\in\Sigma_1$.
On the other hand, we have for any $c>0$ (see [\cite{LiebLoss},\cite{Mossino}]),
$$\int_\R (c-\xi^2)^{+}\mod{\psi(\xi)}^2\,d\xi\le \int_\R (c-\xi^2)^{+}\mod{\psi_*(\xi)^2}\,d\xi,$$
which gives
$$c\int_{-\sqrt{c}}^{\sqrt c}\bigl(\mod{\psi(\xi)}^2-\mod{\psi_*(\xi)}^2\bigr)\,d\xi
   \le \int_{-\sqrt{c}}^{\sqrt c}\xi^2\bigl(\mod{\psi(\xi)}^2-\mod{\psi_*(\xi)}^2\bigr)\,d\xi.\reqlabel{Form*5}$$
Since the symmetric-decreasing rearrangement is order-preserving, it follows that $\psi_*$ also satisfies \cite{DecaiExp}.
Therefore, using the L'Hospital rule we get, for $f(\xi)\colon=\mod{\psi(\xi)}^2-\mod{\psi_*(\xi)}^2$,
$$\lim_{c\to+\infty}c\int_{-\sqrt{c}}^{\sqrt c}f(\xi)\,d\xi=\lim_{c\to+\infty}c^2\bigl[f(-\sqrt{c})-f(\sqrt{c})\bigr]=0$$
and consequently
$$\int_\R\xi^2\mod{\psi_*(\xi)}^2\,d\xi\le \int_\R\xi^2\mod{\psi(\xi)}^2\,d\xi.\reqlabel{Form*6}$$
The first equality of \cite{Form*4} with $p=4$, together with
the second inequality in \cite{Form*4} and \cite{Form*6} imply that $E(\psi_*)=E_\mn$, which means that $\psi_*\in{\cal G}$.
Moreover, for $p=+\infty$, it follows that
$\psi_*(0)=\max\{\mod{\psi(\xi)}\,;\, \xi\in\R\}$, and the proof is complete.\quad\cqd

\smallskip\noindent
{\bf Remark\ \lemlabel{Rem5}:} Theorem~\cite{Thm1} asserts that  ${\cal G}$ is not empty. In fact, ${\cal G}$ has infinitely many elements,
because if $\psi\in{\cal G}$, then $e^{\unim\theta}\psi\in{\cal G}$, for all $\theta\in\R$. Nevertheless, we do not know if $\cal G$ has only
one real positive valued function in the case $\lamb<0$. However, if there are multiple real valued ground states $\psi$ in $\cal G$,
one should be aware that the Lagrange multiplier $\mu$ might depend also on $\psi$. Indeed, multiplying both sides of Eq.~\cite{GPE1} by $\psi$, we get
$$\mu(\psi)=E_\mn+{\lamb \over 2}\int_\R\mod{\psi(\xi)}^4d\xi.\reqlabel{FormMu1}$$

As we are going to see in the study of stability, it is important to notice that the set of Lagrange multipliers is bounded. This is
immediate because
$$\mod{\mu(\psi)}\le \mod{E_\mn}+{\mod{\lamb }\over 2}\norma{\psi}{4}^4\le
  \mod{E_\mn}+{\mod{\lamb} C_{\hbox{\letranota gn}}\over 2}\norma{\psi}{X}$$
and because $\cal G$ is bounded in $X$.

\smallskip\noindent
{\bf Remark\ \lemlabel{Rem7}:}  Theorem~\cite{Thm1} states that Eq.~\cite{GPE1} admits ground state solutions
in both cases: attractive ($\lamb<0$) and repulsive ($\lamb>0$). In the repulsive case, Eq.~\cite{GPE1} presents a different behavior
when compared with the classical NLS. In fact, assume that $\varphi\in H^1(\R)$ is a nontrivial solution of
$$-\varphi''+\lamb\mod{\varphi}^2\varphi=\mu\varphi.\reqlabel{NLS}$$
If we multiply both sides of \cite{NLS} by $\varphi(\xi)$ (respectively $-\xi\varphi'(\xi)$) and integrate on $\R$,
we get
$$\eqalign{
\mu & = \int_\R\mod{\varphi'(\xi)}^2d\xi + \lamb\int_\R\mod{\varphi(\xi)}^4d\xi, \cr
\mu & = -\int_\R\mod{\varphi'(\xi)}^2d\xi + {\lamb\over 2}\int_\R\mod{\varphi(\xi)}^4d\xi. \cr
}$$
Hence,
$$2\int_\R\mod{\varphi'(\xi)}^2d\xi + {\lamb\over 2}\int_\R\mod{\varphi(\xi)}^4d\xi=0,$$
which is impossible if $\lamb\ge0$ and $\varphi\not\equiv 0$. Therefore, \cite{NLS} has no nontrivial solution in $H^1(\R)$ if
$\lamb\ge 0$.

\smallskip
\noindent$\bullet$ $\underline{\hbox{\sl Stability of ground states}}$
\smallskip

In order to prove the stability of ground states of \cite{GPE1}, let us consider the Cauchy problem
$$\unim{\partial v\over \partial\tau}={1\over 2}E'(v),\quad v(0,\xi)=v_0(\xi),\reqlabel{ProbCauchy}$$
where $E'$ is the Fr\'echet derivative of $E$ in $X$, i.e.,
$${1\over 2}E'(v):=-{\partial^2 v\over \partial \xi^2}+\xi^2v
    +\lamb \mod{v}^2v-V_0\cos^2\left(\alpha\xi\right)v.$$

It is well known [\cite{Carles}, \cite{Oh}] that \cite{ProbCauchy} has a unique solution that is global in time, i.e., for any $v_0\in X$, there exists a
unique $v\in C\bigl([0,+\infty),X\bigr)$ satisfying \cite{ProbCauchy}. In particular, if $\psi\in{\cal G}$, the unique solution $u$ of \cite{ProbCauchy} such that
$u(0,\xi)=\psi(\xi)$ is  the standing wave given by
$$u(\tau,\xi)=e^{-\unim\mu\tau}\psi(\xi).$$
Saying that $\psi$ is stable means that if the initial datum $v_0$ of Eq.~\cite{ProbCauchy} is close enough to $\psi$, then the
trajectories $v(\tau,\cdot)$ remain close to the set $\cal G$, as $\tau$ varies in $\R$. More precisely,

\smallskip\noindent
\noindent{\bf Definition\ \lemlabel{DefStab}:} {\sl We will say that $\cal G$ is stable if, for each $\varepsilon>0$, there exists $\delta>0$
such that if $v_0\in X$ satisfies $\inf_{\psi\in\cal G}\norma{v_0-\psi}{X}<\delta$, then the solution of \cite{ProbCauchy} satisfies
$$\sup_{\tau\in\R}\,\inf_{\psi\in{\cal G}}\norma{v(\tau,\cdot)-e^{-\unim\mu\tau}\psi}{X}<\varepsilon.$$  }

Before proving the stability of $\cal G$, we state and prove the following lemma about the compactness  of the set $\cal G$, which will be
needed in the proof of the stability result.

\smallskip\noindent
\noindent{\bf Lemma\ \lemlabel{CompacSetG}:} {\sl The set $\cal G$ is compact and weakly sequentially closed in $X$. More precisely, if
$\{\psi_n\}_{n\in\N}$ is a sequence of $\cal G$, then there exists $\psi\in\cal G$ and a subsequence $\{\psi_{n_k}\}_{k\in\N}$ such that
$\psi_{n_k}\rightarrow\psi$ in $X$. Also, if $\{\psi_n\}_{n\in\N}$ is a sequence of $\cal G$ and $\psi_n\rightharpoonup\psi$ in $X$-weakly,
then $\psi_n\rightarrow\psi$ strongly in $X$.}

\smallskip\noindent{\bf Proof:} Let $\{\psi_n\}_{n\in\N}$ be a sequence of $\cal G$. Since $\cal G$ is bounded in $X$, there exists a
subsequence $\{\psi_{n_k}\}_{k\in\N}$ and $\psi\in X$ such that $\psi_{n_k}\rightharpoonup\psi$ in $X$-weak. By the compactness of the
embedding $X\subset L^p(\R)$ for $2\le p<\infty$, we conclude that $\psi\in\Sigma_1$.
Now, arguing as in the step~2 of the proof of Theorem~\cite{Thm1}, we obtain that $\psi\in\cal G$. To conclude the proof, it remains to show
that $\psi_{n_k}\rightarrow\psi$ in $X$ strongly. To see this, it is enough to remark that $\norma{\psi_{n_k}}{X}\rightarrow\norma{\psi}{X}$,
but this follows from \cite{Part1} and the fact that
$$\norma{\psi_{n_k}}{X}^2=E_\mn-{\lamb \over 2}\int_\R\mod{\psi_{n_k}(\xi)}^4d\xi
     +V_0\int_\R\cos^2\left(\alpha\xi\right)\mod{\psi_{n_k}(\xi)}^2d\xi .\reqlabel{blable}$$
This proves that $\cal G$ is compact and the remaining claims follow easily.\quad\cqd

\smallskip
In order to prove the stability of $\cal G$, we need the following well known conservation laws that hold for all solutions of
\cite{ProbCauchy}: {\sl Assume that $v(\tau,\xi)$ is a solution of \cite{ProbCauchy} such that $v(0,\xi)=v_0(\xi)$. Then, for each $\tau\in\R$, we
have}
$$Q\bigl(v(\tau,\cdot)\bigr) = Q(v_0)\quad\hbox{\rm and}\quad E\bigl(v(\tau,\cdot)\bigr) = E(v_0).$$

Concerning  the above identities, the first one is known as the ``conservation of charge'' and the second one is the ``conservation of
energy''. The first one can be obtained multiplying Eq.~\cite{ProbCauchy} by $-\unim\overline{v(\tau,\xi)}$ (i.e., the complex conjugate of
$\unim v(\tau,\xi)$). By the same way, one obtains the conservation of energy, but in that case we should multiply the equation by
${\partial\overline{v}\over\partial\tau }$.

We are now in position to prove the stability of $\cal G$.

\smallskip\goodbreak
\noindent{\bf Theorem\ \lemlabel{Thm2}:} {\sl The set $\cal G$ is stable in the sense of the Definition~\cite{DefStab} .}

\smallskip\noindent
{\bf Proof:}  Arguing by contradiction, we proceed as in Cazenave-Lions [\cite{CazeLions}]. If $\cal G$ is not stable, there exists $\varepsilon_0>0$ such
that for all integer $n\in\N$, we can find  $v_{0n}\in X$ satisfying
$$r_n:=\inf_{\psi\in\cal G}\norma{v_{0n}-\psi}{X}<{1\over n}\reqlabel{PrimaCond}$$
and
$$\sup_\tau\inf_{\psi\in\cal G}\norma{v_n(\tau ,\cdot)-e^{-\unim\mu(\psi)\tau}\psi}{X}\ge \varepsilon_0,\reqlabel{DuaCond}$$
where $v_n\in C(\R;X)$ is the unique solution of \cite{ProbCauchy} with initial datum $v_{0n}$.

Let $\psi_n\in\cal G$ such that $r_n\le \norma{v_{0n}-\psi_n}{X}<r_n+1/n$. Since $\cal G$ is bounded in $X$ and the sequence $\mu(\psi_n)$ is
bounded in $\R$ (see Remark~\cite{Rem5}),  there exists $(\psi_\infty,\mu_\infty)\in X\times\R$ and a subsequence still denoted by
$(\psi_{n},\mu(\psi_{n}))$ such that $\psi_{n}$ converges to $\psi_\infty$ weakly in $X$ and strongly in $L^4(\R)\cap L^2(\R)$, while
$\mu(\psi_{n})$ converges to $\mu_\infty$ in $\R$. By Lemma~\cite{CompacSetG} we know that $\psi_{n}\rightarrow\psi_\infty$ in $X$ and $\psi_\infty\in\cal G$.
Also, from the fact that $\norma{v_{0n}-\psi_n}{X}\rightarrow 0$, we may infer that $v_{0n}$ converges to $\psi_\infty$ in $X$ and in $L^4(\R)\cap
L^2(\R)$ as well.

On the other hand, one may observe from \cite{DuaCond} that there exists $\tau_n\in\R$ such that
$$\inf_{\psi\in\cal G}\norma{v_n(\tau_n,\cdot)-e^{-\unim\mu(\psi)\tau_n}\psi}{X}\ge{1\over 2}\varepsilon_0.\reqlabel{ParaContrad}$$

Setting $\widetilde{\psi}_n:=e^{\unim\mu(\psi_n)\tau_n}v_{n}(\tau_{n},\cdot)$, it follows from the conservation of charge and energy that
$$Q(\widetilde\psi_{n})=Q(v_{0n})\rightarrow Q(\psi_\infty)=1,\quad E(\widetilde\psi_n)=E(v_{0n})\rightarrow E(\psi_\infty)=E_\mn
     \quad\hbox{\rm as}\quad n\rightarrow\infty.$$

All this means that $\{\widetilde\psi_n\}_{n\in\N}$ is a bounded sequence in $X$ and there exists $\widetilde\psi_\infty\in X$ and a subsequence
(still denoted by $\{\widetilde\psi_n\}_{n\in\N}$) converging to $\widetilde\psi_\infty$ weakly in $X$ and strongly in $L^4(\R)\cap L^2(\R)$.
Therefore, $\widetilde\psi_\infty\in\cal G$ and $E(\widetilde\psi_n)\rightarrow E(\widetilde\psi_\infty)$. Again, invoking relation \cite{blable} we observe
that $\widetilde\psi_n\rightarrow\widetilde\psi_\infty$ strongly in $X$ and this is a contradiction with \cite{ParaContrad}. This finishes the proof.\quad\cqd

\smallskip
\noindent$\bullet$ $\underline{\hbox{\sl The minimal energy as function of $\lamb$}}$
\smallskip

In order to explicit the dependence of the minimal energy relatively to the parameter $\lamb$, let us denote
$$E_\mn(\lamb):=\min\Bigl\{E_\lamb(\psi)\,;\, \psi\in\Sigma_1\Bigr\},\reqlabel{2aa}$$
where $E_\lamb$ is the energy functional introduced in \cite{DefEQ}, and
${\cal G}_\lamb:=\Bigl\{\psi\in\Sigma_1\,;\, E_\lamb(\psi)=E_\mn(\lamb)\Bigr\}$
the set of corresponding ground-states. Theorem~\cite{Thm1} assures that $E_\mn$ is well-defined as function of $\lamb\in\R$
and it is easy to see that it is strictly increasing. Indeed, for $h>0$ and $\psi\in{\cal G}_{\lamb+h}$ we have
$$E_\mn(\lamb)\le E_\lamb(\psi)=E_\mn(\lamb+h)-{h\over 2}\norma{\psi}{4}^4<E_\mn(\lamb+h).\reqlabel{Cresc}$$

\smallskip\goodbreak
\noindent{\bf Proposition\ \lemlabel{Prop}:} {\sl $E_\mn(\lamb)$ is a strictly increasing and concave function such that
$$\lim_{\lamb\to\pm\infty}E_\mn(\lamb)=\pm\infty.\reqlabel{LimInfinito}$$
}

The proof of the Proposition relies on the following:
\goodbreak
\noindent{\bf Lemma\ \lemlabel{UnifLimit}:} {\sl For each $a,b\in\R$, $a<b$, the set $\bigcup_{a\le\lamb\le b}{\cal G}_\lamb$ is bounded in $X$.
More precisely, there exists a constant $C_{a,b}>0$ such that
$$\norma{\psi}{X}\le C_{a,b}\quad\forall \psi\in\bigcup_{a\le\lamb\le b}{\cal G}_\lamb.$$}

\noindent
{\bf Proof:} It suffices to prove for $a<0<b$. Let $\lamb\in[a,b]$ and $\psi\in{\cal G}_\lamb$. Then
$$E_\mn(b)\ge E_\mn(\lamb)=E_\lamb(\psi)\ge \norma{\psi}{X}^2+{\lamb\over 2}\norma{\psi}{4}^4-\mod{V_0}.$$
If $\lamb\ge 0$, then $\norma{\psi}{X}^2\le E_\mn(b)+\mod{V_0}$.
If $\lamb<0$, it follows from Gagliardo-Nirenberg inequality \cite{Gal-Nir} that
$$E_\mn(b)\ge \norma{\psi}{X}^2+{\lamb C_{\hbox{\letranota gn}}\over 2}\norma{\psi}{X}-\mod{V_0}\ge \norma{\psi}{X}^2+{aC_{\hbox{\letranota gn}}\over 2}\norma{\psi}{X}-\mod{V_0}$$
and we have
$$\norma{\psi}{X}^2\le 2\left(E_\mn(b)+{a^2C_{\hbox{\letranota gn}}^2\over 8}+\mod{V_0}\right).\quad\cqd$$
\smallskip\noindent
{\bf Proof of Proposition~\cite{Prop}:} It follows from \cite{Cresc} that
$$E_\mn(\lamb)\le \liminf_{h\to 0^+}E_\mn(\lamb+h).\reqlabel{2bb}$$
Let $\psi\in{\cal G}_\lamb$. Then,
$$E_\mn(\lamb+h)\le E_{\lamb+h}(\psi)=E_\mn(\lamb)+{h\over 2}\norma{\psi}{4}^4\,\,\,\forall h\in\R,\reqlabel{2bc}$$
from which we obtain
$$\limsup_{h\to 0}E_\mn(\lamb+h)\le E_\mn(\lamb).\reqlabel{2cc}$$
From \cite{2bb} and \cite{2cc} we conclude that
$$\lim_{h\to 0^+}E_\mn(\lamb+h)=E_\mn(\lamb).\reqlabel{2ccc}$$
On the other hand, if $h<0$ and $\psi\in{\cal G}_{\lamb+h}$, we have
$$E_\mn(\lamb)\le E_\lamb(\psi)=E_\mn(\lamb+h)-{h\over 2}\norma{\psi}{4}^4,\reqlabel{cd}$$
from which, using the Gagliardo-Nirenberg inequality we obtain,
$$E_\mn(\lamb)\le E_\mn(\lamb+h)-{hC_{\hbox{\letranota gn}}\over 2}\norma{\psi}{X}.\reqlabel{cdd}$$
By choosing and fixing $a\in\R$ such that $a<\lamb+h<\lamb$, it follows from Lemma~\cite{UnifLimit}	that there exists a constant $C_{a,\lamb}$
such that
$$E_\mn(\lamb)\le E_\mn(\lamb+h)-{hC_{\hbox{\letranota gn}}\over 2}C_{a,\lamb},$$
from which we get
$$\liminf_{h\to 0^-}E_\mn(\lamb+h)\ge E_\mn(\lamb).\reqlabel{2dd}$$
So, \cite{2cc} and \cite{2dd} give
$$\lim_{h\to 0^-}E_\mn(\lamb+h)=E_\mn(\lamb)$$
and we conclude from \cite{2ccc} that $E_\mn$ is a continuous function. Moreover, from \cite{2bc} we have
$$E_\mn(\lamb+h)+E_\mn(\lamb-h)-2E_\mn(\lamb)\le 0\,\,\,\forall\lamb,h\in\R,$$
which, under continuity, is sufficient to assure the concavity of $E_\mn$.

In order to prove \cite{LimInfinito}, consider $\psi_k\in{\cal G}_k$, $k\in\N$. Then, using \cite{Cresc} and \cite{2bc} with $h=1$
we have
$${1\over 2}\norma{\psi_{k+1}}{4}^4 \le E_\mn(k+1)-E_\mn(k) \le {1\over 2}\norma{\psi_k}{4}^4,\quad\forall k\in\N$$
from which we obtain for all $n\in\N$
$${1\over 2}\sum_{k=1}^{n+1}\norma{\psi_k}{4}^4 \le E_\mn(n+1)-E_\mn(0) \le {1\over 2}\sum_{k=0}^n\norma{\psi_k}{4}^4.$$
Now, arguing by contradiction, assume that there exists a positive constant $C$ such that $E_\mn(\lamb)\le C$ for all $\lamb\in\R$. Then,
it follows from the first inequality in the above expression that the sequence $\{\psi_n\}_{n}$ converges to zero in $L^4(\R)$.
But, repeating the arguments used in the first step of the proof of Theorem~\cite{Thm1}, we conclude that the sequence $\{\psi_k\}_k$ is bounded in $X$,
so that, passing to a subsequence if necessary, we have that $\psi_k$ converges to $\psi$ weakly in $X$. Since the embedding $X\subset L^p(\R)$
is compact for $2\le p<\infty$, we conclude, by choosing $p=4$ that $\psi=0$,  and by choosing $p=2$ that $\psi\in\Sigma_1$, which is absurd.
Since the same arguments apply for $\psi_{-k}\in{\cal G}_{-k}$, $k\in\N$, we finish the proof.\quad\cqd

\smallskip\noindent
{\bf Remark\ \lemlabel{RemMuLamb}:} Concerning the dependence of the chemical potential relatively to $\lamb$, it follows from the formula \cite{FormMu1} and the
characterization of ${\cal G}_\lamb$ given by Theorem~\cite{Thm1-1} that we can also define, for $\lamb\ge0$, the function $\mu_\mn(\lamb)$.
However, as reported in Remark~\cite{Rem5}, it is unclear that $\mu_\mn(\lamb)$ be uniquely determined if $\lamb<0$.
No matter whether or not it be uniquely determined, the fact that
$\lamb\bigl(\mu_\lamb(\psi)-E_\mn(\lamb)\bigr)\ge 0$ for all $\lambda\in\R$ and for all $\psi\in{\cal G}_\lamb$, allows  to state that
$$\lim_{\lamb\to\pm\infty}\mu_\lamb(\psi_\lamb)=\pm\infty,\quad \forall\psi_\lamb\in{\cal G}_\lamb.$$

\smallskip
\noindent{\bf Corollary\ \lemlabel{DecresNorma4}:} {\sl Let $\psi_i\in{\cal G}_{\lamb_i}$, $i=1,2$. If $\lamb_1<\lamb_2$, then $\norma{\psi_2}{4}\le\norma{\psi_1}{4}$}.

\smallskip\noindent
{\bf Proof:} Let $h=\lamb_2-\lamb_1$. Then, from \cite{Cresc} and \cite{2bc}, it follows that
$${1\over 2}\norma{\psi_2}{4}^4 \le {E_\mn(\lamb_2)-E_\mn(\lamb_1)\over  h} \le {1\over 2}\norma{\psi_1}{4}^4.\quad\cqd$$


\bigskip
\noindent{\bf\numsection.\ A new method to obtain approximations for the minimal energy}\par
\bigskip

The results presented in the previous section, although mathematically rigorous, do not
provide sufficiently precise quantitative information for some relevant physical quantities.
Since it seems to be not possible to calculate exact explicit solutions of Eq.~\cite{GPE1} (in particular,
the ground state solution $\psi_\mn$), the exact values of such quantities cannot be expressed
in terms of the known parameters. Therefore, it would be useful to obtain some kind of explicit formulae through which
one could approximate them. Indeed, these are outcomes mainly interesting from the physical point of view.
In this sense, there exist some well known methods that were already applied to the problem we are dealing with
(see [\cite{TrallB},\cite{TrallD},\cite{TrallA}]). Nevertheless, as we shall see below, a new interesting and powerful
approach can be developed.

Let $\{\Fi_\lamb\}_{\lamb\in\R}$ a family of functions in $\Sigma_1$ such that $\lamb\mapsto\Fi_\lamb$ defines a differentiable curve in $X$.
Then,
$$\int_\R\Fi_\lamb(\xi){d\hfil\over d\lamb}\Fi_\lamb(\xi)\,d\xi=0\,\quad\forall\lamb\in\R.\reqlabel{Ort}$$
Since the energy functional $E_\lamb$ is differentiable in $X$, it follows from the chain rule that $\lamb\mapsto E_\lamb(\Fi_\lamb)$
is also a differentiable function and
$${d\hfil\over d\lamb}E_\lamb(\Fi_\lamb) =\Bigl\langle E_\lamb'(\Fi_\lamb)\,|\,{d\hfil\over d\lamb}\Fi_\lamb\Bigr\rangle
               +{1\over 2}\int_\R\mod{\Fi_\lamb(\xi)}^4d\xi,\reqlabel{Form3.2}$$
where $\langle\cdot |\,\cdot\rangle$ denotes the duality product between $X$ and its dual $X^*$.

From now on we assume that, for any $\lamb\in\R$, one can choose a real ground state $\psi_\lamb\in {\cal G}_\lamb$ such that $\lamb\mapsto\psi_\lamb$ defines
a differentiable curve in $X$. In this case,
$$E_\mn(\lamb)\colon=\int_\R\mod{\psi_\lamb'(\xi)}^2d\xi+\int_\R\xi^2\mod{\psi_\lamb(\xi)}^2d\xi
              +{\lamb\over 2}\int_\R\mod{\psi_\lamb(\xi)}^4d\xi-V_0\int_\R\cos^2(\alpha\xi)\mod{\psi_\lamb(\xi)}^2d\xi $$
is differentiable as a function of $\lamb$ and, as $E_\lamb'(\psi_\lamb)=\mu\psi_\lamb$, we get from \cite{Ort} and
\cite{Form3.2}
$${d\hfil\over d\lamb}E_\mn(\lamb)={1\over 2}\int_\R\mod{\psi_\lamb(\xi)}^4d\xi={1\over 2}\norma{\psi_\lamb}{4}^4\reqlabel{DerivEnerg} $$
and we have the formula
$$E_\mn(\lamb)=E_\mn(0)+{1\over 2}\int_0^\lamb\norma{\psi_s}{4}^4\,ds.  \reqlabel{Form3.3}$$
On the other hand, if we denote $\mu_\mn(\lamb)\colon= \mu(\psi_\lamb)$,  it follows from \cite{FormMu1} that
$${d\hfil\over d\lamb}\mu_\mn(\lamb) =\norma{\psi_\lamb}{4}^4+{\lamb\over 2}{d\hfil\over d\lamb}\norma{\psi_\lamb}{4}^4,\reqlabel{Form3.4}$$
and we get by integration on $\lamb$ the formula
$$\mu_\mn(\lamb)=\mu_\mn(0)+{1\over 2}\left(\lamb\norma{\psi_\lamb}{4}^4+\int_0^\lamb\norma{\psi_s}{4}^4\,ds\right).\reqlabel{Form3.5}$$

Notice that $E_\mn(0)=\mu_\mn(0)$ and, in the case $V_0=0$, we have from Remark \cite{Rem3} that $E_\mn(0)=\mu_\mn(0)=1$.

\smallskip\goodbreak
\noindent{\bf Remark\ \lemlabel{NewFormula}:} Formul\ae\ \cite{Form3.3} and \cite{Form3.5} express the minimal energy $E_\mn$
and the corresponding chemical potential $\mu_\mn$ as functions of $\lamb$, which also depend explicitly on the $L^4$-norm
of unknown ground states. However, by eliminating this explicit dependence,
we can obtain an exact formula relating these two quantities. More precisely, by using \cite{DerivEnerg}, we can rewrite \cite{Form3.4} as
$${d\hfil\over d\lamb}\mu_\mn(\lamb)=2{d\hfil\over d\lamb}E_\mn(\lamb)+\lamb{d^2\hfil\over d\lamb^2}E_\mn(\lamb),$$
from which we get easily
$${d\hfil\over d\lamb}\left(E_\mn(\lamb)-{1\over\lamb}\int_0^\lamb\mu_\mn(s)\,ds\right)=0.$$
Hence, there exists a constant $C$ such that
$$E_\mn(\lamb)-{1\over\lamb}\int_0^{\lamb}\mu_\mn(s)\,ds=C,\,\,\,\forall\lamb\in\R.$$
Since $E_\mn(0)=\mu_\mn(0)$, it follows that $C=0$ and we have the identity
$$E_\mn(\lamb)={1\over\lamb}\int_0^{\lamb}\mu_\mn(s)\,ds,\,\,\,\forall\lamb\in\R.\reqlabel{NewFormEMu}$$

\smallskip
The formul\ae\ \cite{Form3.3} and \cite{Form3.5} can be used to obtain explicit approximate functions depending on $\lamb$
for the minimal energy $E_\mn$ and the corresponding chemical potential $\mu_\mn$, by choosing appropriate trial
functions that generate differentiable curves in $\Sigma_1$.

Motivated by inequality \cite{DecaiExp}, we consider the trial functions $\phi_\kappa\in\Sigma_1$ defined as
$$\phi_\kappa(\xi)\colon=\root 4 \of {2\kappa\over\pi}\exp(-\kappa\xi^2).$$
A direct calculation shows that the function $\kappa\mapsto E_\lamb(\phi_\kappa)$ is given by
$$E_\lamb(\phi_\kappa)=\kappa+{1\over 4\kappa}+{\lamb\sqrt{\kappa}\over 2\sqrt\pi}
  -{V_0\over 2}\left(1-\exp\left(-{\alpha^2\over 2\kappa}\right)\right).$$
We can show that, for all $\lamb\in\R$, the value of $\kappa$ that minimizes the function
$\kappa\mapsto E_\lamb(\phi_\kappa)$ is given by the largest solution (in fact, unique solution if $\alpha^2V_0<1$)
$\kappa(\lamb)$ of the following transcendental equation
$$\kappa^{3/2}\left(\kappa^{1/2}+{\lamb\over 4\sqrt\pi}\right)
    +{V_0\alpha^2\over 4}\exp\left({-\alpha^2 \over 2\kappa}\right)={1\over 4}.\reqlabel{EqAlg1}$$

In order to obtain explicit approximate formul\ae\ for the minimal energy $E_\mn(\lamb)$ and the corresponding
chemical potential $\mu_\mn(\lamb)$, we introduce the functions
$$\left\{\eqalign{
E_{\hbox{\fiverm app}}(\lamb)   & \colon=E_\mn(0)+{1\over 2}\int_0^\lamb\norma{\Fi_s}{4}^4ds,\cr
\mu_{\hbox{\fiverm app}}(\lamb) & \colon=\mu_\mn(0)+{1\over 2}\left(\lamb\norma{\Fi_\lamb}{4}^4+\int_0^\lamb\norma{\Fi_s}{4}^4\,ds\right),\cr
}\right.$$
where, in this case,
$$\Fi_\lamb(\xi)\colon=\phi_{\kappa(\lamb)}(\xi)=\left({2\kappa(\lamb)\over \pi}\right)^{1/4}\exp(-\kappa(\lamb)\xi^2).$$
By a straightforward calculation we get
$$\left\{\eqalign{
E_{\hbox{\fiverm app}}(\lamb) & = E_\mn(0)+{1\over 2\sqrt{\pi}}\int_0^\lamb\sqrt{\kappa(s)}\,ds, \cr
\mu_{\hbox{\fiverm app}}(\lamb) & = \mu_\mn(0)+{1\over 2\sqrt{\pi}}\left(\lamb\sqrt{\kappa(\lamb)}+\int_0^\lamb\sqrt{\kappa(s)}\,ds\right).\cr
}\right.\reqlabel{FormAprox}$$
Moreover, by arguing as in Remark~\cite{NewFormula}, we can show that
$$E_{\hbox{\fiverm app}}(\lamb)={1\over\lamb}\int_0^\lamb \mu_{\hbox{\fiverm app}}(s)\,ds.$$
From the above identity we can easily relate the Taylor coefficients $E_n$ of $E_{\hbox{\fiverm app}}$ with the ones of
$\mu_{\hbox{\fiverm app}}$. In fact, we have $\mu_n=(n+1)E_n$ for all $n\in\N$.

A relatively simple situation is the one at which the optical lattice potential is not present ($V_0=0$).
In this case, \cite{EqAlg1} has a unique solution and if we define $\sigma(\lamb)\colon=\sqrt{\kappa(\lamb)}$,
the equation \cite{EqAlg1} with $V_0=0$ can be written as
$$\sigma^4+{\lamb\over 4\sqrt\pi}\sigma^3={1\over 4},\,\,\, \forall\lamb\in\R.\reqlabel{NovaVariavel}$$
By differentiating this equation implicitly with respect to $\lamb$,
we get easily the following properties of the function $\sigma(\lamb)$:

\smallskip\goodbreak
\noindent{\bf Lemma\ \lemlabel{PropFuncSigma}:} {\sl The function $\sigma:\R\rightarrow\R$ is $C^\infty$, positive,
strictly decreasing, convex and satisfies the following properties:
$$\sigma(0)=\sqrt2/2;\quad \lim_{\lamb\to-\infty}\sigma(\lamb)=+\infty\quad\hbox{\sl and}\quad\lim_{\lamb\to+\infty}\sigma(\lamb)=0.$$
More precisely,
$$\sigma(\lamb)\sim -{\lamb\over 4\sqrt\pi}\,\,\,\hbox{\sl as}\,\,\,\lamb\rightarrow-\infty\quad\hbox{\sl and}\quad
   \sigma(\lamb)\sim\left(\sqrt\pi\over\lamb\right)^{1/3}\,\,\,\hbox{\sl as}\,\,\,\lamb\rightarrow+\infty.\quad\cqd$$
}

As consequence of the previous lemma, we can show that $E_{\hbox{\fiverm app}}(\lamb)$ satisfies the
general properties of $E_\mn(\lamb)$ as reported in Proposition~\cite{Prop}. More precisely,

\smallskip\goodbreak
\noindent{\bf Corollary\ \lemlabel{PropFuncEMu}:} {\sl The functions
$E_{\hbox{\fiverm app}}(\lamb)$ and $\mu_{\hbox{\fiverm app}}(\lamb)$ defined in \cite{FormAprox}
are $C^\infty$, strictly increasing and concave. Moreover
$$\left\{\eqalign{
  & E_{\hbox{\fiverm app}}(\lamb)\,\, \hbox{\sl and}\,\,\,\mu_{\hbox{\fiverm app}}(\lamb)\,\,\,\hbox{\sl are}\,\,\,O(-\lamb^{2})\,\,\hbox{\sl as}\,\,\lamb\rightarrow-\infty,\cr
  & E_{\hbox{\fiverm app}}(\lamb)\,\, \hbox{\sl and}\,\,\,\mu_{\hbox{\fiverm app}}(\lamb)\,\,\,\hbox{\sl are}\,\,\,O(\lamb^{2/3})\,\,\hbox{\sl as}\,\,\lamb\rightarrow+\infty.\cr
  }\right.\quad\cqd$$

}

By implicit differentiation the Eq.~\cite{NovaVariavel} on $\lamb$, we get
$$\sigma(0)={\sqrt{2}\over 2},\,\,\,\sigma'(0)=-{1\over 16\sqrt\pi},\,\,\,\sigma''(0)={3\over 128\pi\sqrt2},\,\,\,
\sigma'''(0)=-{3\over 512\pi\sqrt\pi},\,\,\, \sigma^{(4)}(0)={45\over 16384\pi^2\sqrt2},\cdots\reqlabel{Derivadas}$$    
Then, for $\lamb$ small enough, we can consider the approximation
$$\sigma(\lamb)\approx{\sqrt{2}\over 2}-{\lamb\over 16\sqrt\pi}+{3\lamb^2\over 256\pi\sqrt2}
-{\lamb^3\over 1024\pi\sqrt\pi}+{45\lamb^4\over 393216\pi^2\sqrt2}.$$
By substituting the above approximation in \cite{FormAprox} we have, for $\lamb$ small enough,
$$\left\{\eqalign{
E_{\hbox{\fiverm app}}(\lamb) & \approx 1+{\lamb\over 2\sqrt{2\pi}}-{\lamb^2\over 64\pi}+{\lamb^3\over 512\pi\sqrt{2\pi}}
   -{\lamb^4\over 8192\pi^2}+{9\lamb^5\over 786432\pi^2\sqrt{2\pi}},\cr
\mu_{\hbox{\fiverm app}}(\lamb) & \approx 1+{\lamb\over \sqrt{2\pi}}-{3\lamb^2\over 64\pi}+{4\lamb^3\over 512\pi\sqrt{2\pi}}
   -{5\lamb^4\over 8192\pi^2}+{54\lamb^5\over 786432\pi^2\sqrt{2\pi}}. \cr
}\right.\reqlabel{FormAproExpl}$$

\smallskip
\noindent{\bf Remark\ \lemlabel{RemTrall}:}
The above results are valid for attractive ($\lamb <0$) as well as for repulsive ($\lamb >0$) interatomic interaction strengths.
By using \cite{DerivEnerg}, \cite{Form3.4} and \cite{PsiZero} we can show that $E'_\mn(0)=E'_{\hbox{\fiverm app}}(0)=1/2\sqrt{2\pi}$,
which implies that the formul\ae\ given by \cite{FormAproExpl} coincide with $E_\mn(\lamb)$ and $\mu_\mn(\lamb)$, respectively,
up to first order terms in $\lamb$.
Notice that, up to second order terms in $\lamb$, the approximate chemical potential $\mu_{\hbox{\fiverm app}}$ can be written as
$$\mu_{\hbox{\fiverm app}}(\lamb)\approx 1+{\lamb\over\sqrt{2\pi}}-\eps\lamb^2,\reqlabel{MuMe}$$
where $\eps={3\over 64\pi}\approx 0.0149207$.

It is noteworthy to compare Eq.~\cite{MuMe} with the one obtained in [\cite{TrallB}]
by a perturbative method. Considering that our dimensionless parameters are in fact twice the ones used there, both formul\ae\
\cite{MuMe} and Eq.~(31) in [\cite{TrallB}] coincide, except for the values of $\eps$: $\eps=0.0149207$ and $\eps=0.016553$ respectively,
which are also very close to each other.
Moreover, Fig.~1 displays a comparison between the values of $\mu_{\hbox{\fiverm app}}$ from Eq.~\cite{FormAproExpl} (solid line) and those obtained
by perturbation theory (dashed line) [\cite{TrallB}]. Also and for sake of comparison, the numerical evaluation of Eq.~\cite{GPE1} is shown
by full stars. From the figure is observed that Eq.~\cite{FormAproExpl} increases the accuracy of the solution with respect to the
perturbative method. In the scale of the figure no significant differences are observed between the numerical solutions and
the calculated values using Eq.~\cite{FormAproExpl} for the interval range $\mod{\lamb}<8$.    
Nevertheless the perturbation method gives a large error for $\mod{\lamb}>6$.               

\setbox11=\hbox{\epsfbox{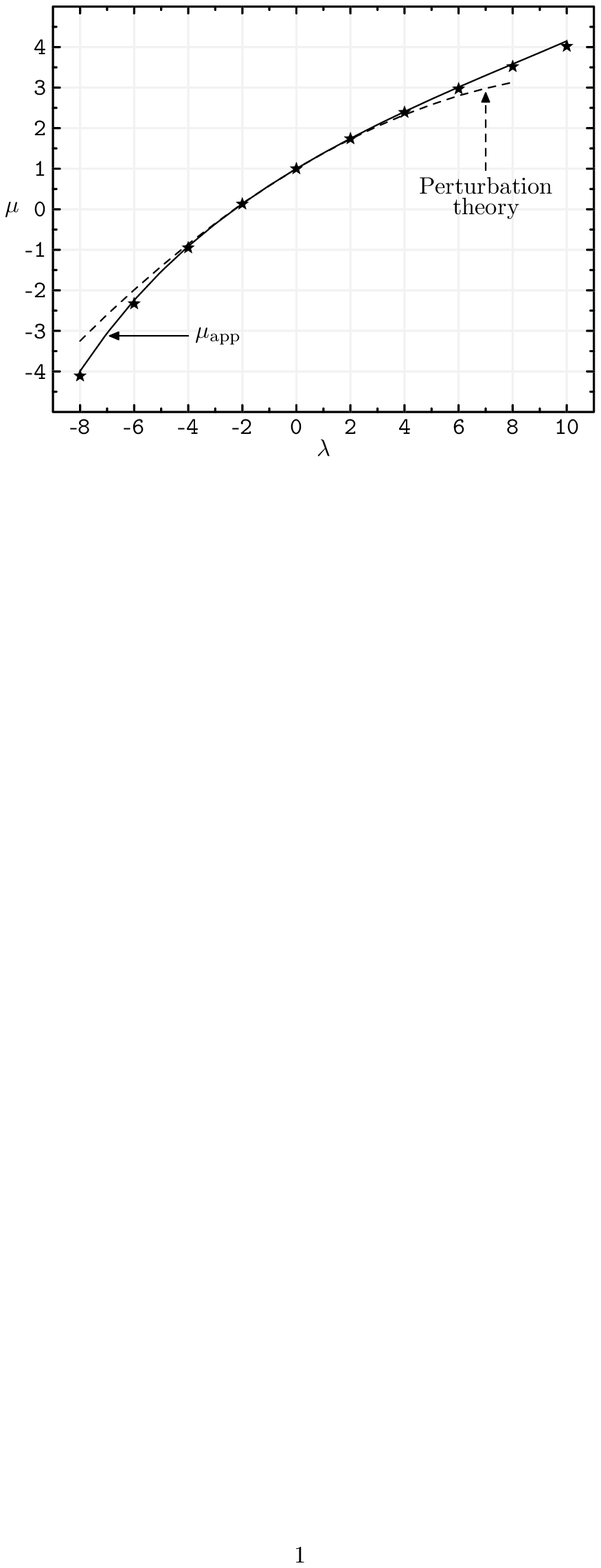}}
$$\hbox{\box11}$$
{\parindent=1cm\narrower\noindent{\bf Fig.1:} \letranota Dimensionless chemical potential $\mu =2\mu _{0}/\hbar \omega$ as a
function of $\lambda=2\lambda _{1D}/(l\hbar \omega )$. Solid line: calculation following Eq.~\cite{FormAproExpl}.
Dashed line: perturbation theory from Ref.~[\cite{TrallB}]. Stars: numerical evaluation of Eq.~\cite{GPE1}.\par}

\noindent
Furthermore, in [\cite{TrallD}] a closed expression for the order parameter is given by (in the case $V_0=0$)
$$\psi(\xi)={1\over\root 4 \of \pi}\exp(-{\xi^2\over 2})\left[1+
   {\lamb\over\sqrt{2\pi}}\int_1^{\sqrt{2}/2}
   \left[{\exp\left(-\xi^2(1-z^2)z^{-2}\right)-z\over 1-z^2}\right]\,dz\right].$$
In our approach, we propose the function of $\Sigma_1$:
$$\Fi_\lamb(\xi)=\left({2\kappa(\lamb)\over\pi}\right)^{1/4}\exp(-\kappa(\lamb)\xi^2),\reqlabel{ApproxGS}$$
where $\kappa(\lamb)$ is the unique root of the  equation \cite{EqAlg1} with $V_0=0$. Recalling that
$\kappa(\lamb)=\sigma(\lamb)^2$, it follows from \cite{Derivadas} that, for $\lamb$ small enough,
$$\left\{\eqalign{
   \kappa(\lamb) & \approx{1\over 2}-{1\over 16}\sqrt{2\over\pi}\lamb+{1\over 128\pi}\lamb^2,\cr
   \root 4\of {\kappa(\lamb)} & \approx {1\over\root 4 \of 2}-{\root 4 \of 2\over 32\sqrt\pi}\lamb+{\root 4\of 2\over 1024\pi\sqrt2}\lamb^2\cr
}\right.$$
and \cite{ApproxGS} can be approximated up to second order terms in $\lamb$ by
$$\widetilde\Fi_{\hbox{\fiverm app}}(\xi)={1\over \root 4\of\pi}\exp\left(-{\xi^2\over 2}\right)\left[1-{\sqrt 2\lamb\over 32\sqrt\pi}+
    {\lamb^2\over 1024\pi}\right]\exp\left[\left({\lamb\over 8\sqrt{2\pi}}-{\lamb^2\over 128\pi}\right)\xi^2\right].\reqlabel{ApproxGSTrunc}$$

It is also interesting to notice that, for all $\lamb\in\R$,
$$\left\{\eqalign{
 {1\over 2}-{\lamb\over 8\sqrt{2\pi}}+{\lamb^2\over 128\pi}& \ge {1\over 4},\cr
 1-{\sqrt 2\lamb\over 32\sqrt\pi}+{\lamb^2\over 1024\pi} & \ge {1\over 2},\cr
}\right.$$
which implies that $\widetilde\Fi_{\hbox{\fiverm app}}$ is a positive function of $X$ for any $\lamb\in\R$.

\smallskip
\noindent{\bf Remark\ \lemlabel{RemStabil}:} It follows from the properties stated in Lemma~\cite{PropFuncSigma} and the Labesgue Theorem
that
$$\lim_{\lamb\to 0}\norma{\Fi_\lamb-\Fi_0}{X}=\lim_{\lamb\to 0}\norma{\widetilde\Fi_\lamb-\Fi_0}{X}=0,$$
where (always assuming that $V_0=0$) $\Fi_0$, $\Fi_\lamb$ and $\widetilde\Fi_\lamb$ are given by \cite{PsiZero}, \cite{ApproxGS} and
\cite{ApproxGSTrunc}, respectively. Therefore, from Theorem~\cite{Thm2}, for $\lamb$ small enough and up to a change of phase, the unique
solutions of Eq.~\cite{EqEvolution} with initial data $\Fi_\lamb$ and $\widetilde\Fi_\lamb$ respectively, remain close (in the sense Definition~\cite{DefStab})
to $u(\tau,\xi):=e^{-\unim\tau}\Fi_0(\xi)$ in $X$, for all time $\tau\in\R$.


\bigskip
\noindent{\bf\numsection.\ Conclusions}\par
\bigskip

Our main results in the first part of this work concern qualitative properties of the minimal energy solutions of 1D GP equation
with cubic nonlinearity in a harmonically
confined periodical potential. We prove the existence of ground states for any $\lamb$ and $V_0$ (Theorem~\cite{Thm1}).
Regardless the value of the laser intensity $V_0$, such ground states have a
Gaussian-like exponential asymptotic behavior, as was pointed out in Theorem \cite{ExpDecay}. For $\lamb>0$ (repulsive interatomic forces),
we prove that there exists a unique positive and symmetric ground state $\psi_\mn$, which is decreasing for $\xi>0$ if $V_0=0$,
and that any other solution of the problem differs from it in a phase factor (Theorem~\cite{Thm1-1}).  We also prove that, independently of the
value of $V_{0}$ and for any $\lamb\in\R$, the set ${\cal G}_\lamb$ of ground states is orbitally stable (Theorem~\cite{Thm2}).
An important consequence of this fact is the stability of those physical magnitudes, which are described by operators defined in the 
Hilbert space $X$ (see \cite{DefX}). This result is particularly related to the superfluidity properties, among other physical phenomena, 
of the harmonically confined condensates loaded in optical lattices [\cite{Burger}].

In the second part (Section~3) we present a new and simple method to construct formul\ae\ (see \cite{Form3.3} and \cite{Form3.5}) 
that allows to approximate the minimal energy, the corresponding chemical potential as well as the ground state. The functions described 
by these formul\ae\ (see \cite{FormAprox}) preserve some global properties of
$E_\mn(\lamb)$ and $\mu_\mn(\lamb)$, as pointed out by Proposition~\cite{Prop} and Corollary~\cite{PropFuncEMu}. In the case $V_0=0$,
we obtain approximations for $E_\mn(\lamb)$ and $\mu_\mn(\lamb)$ as Taylor polynomials of $E_{\hbox{\fiverm app}}(\lamb)$ and 
$\mu_{\hbox{\fiverm app}}(\lamb)$ respectively (see \cite{FormAproExpl}).

\bigskip\goodbreak
\noindent{\bf References}\par
\kern-.6truecm
\MakeBibliography{}

\bye